\documentclass[11pt, reqno]{amsart}
\usepackage{amsmath,amsthm,amssymb,graphicx}

\theoremstyle{plain}
\newtheorem{thm}{Theorem}[section]
\newtheorem{pro}[thm]{Proposition}

\newtheorem{cor}[thm]{Corollary}
\newtheorem{tm}[thm]{Theorem}

\theoremstyle{definition}
\newtheorem{rem}[thm]{Remark}

\numberwithin{equation}{section}

\newcommand{\C}{\mathbb{C}}
\newcommand{\R}{\mathbb{R}}

\newcommand{\cg}{\mathcal{G}}
\newcommand{\cu}{\mathcal{U}}
\newcommand{\ck}{\mathcal{K}}
\newcommand{\cp}{\mathcal{P}}
\newcommand{\co}{\mathcal{O}}

\newcommand{\cs}{\mathcal{S}}

\newcommand{\rmi}{\mathrm{ \, i \, }}
\newcommand{\ud}{\mathrm{ \, d}}
\newcommand{\pa}{\partial}

\newcommand{\n}{\ \vert\ }
\newcommand{\I}{\mathrm{ I \, } }
\newcommand{\II}{\mathrm{ I \! I \, } }

\DeclareMathOperator{\GL}{\mathrm{GL}}

\DeclareMathOperator{\U}{\mathrm{U}}

\DeclareMathOperator{\OO}{\mathrm{O}}

\DeclareMathOperator{\gl}{\mathrm{gl}}

\DeclareMathOperator{\uu}{\mathrm{u}}

\DeclareMathOperator{\oo}{\mathrm{o}}

\DeclareMathOperator{\re}{\mathrm{Re}}
\DeclareMathOperator{\im}{\mathrm{Im}}

\DeclareMathOperator{\diag}{\mathrm{diag}}

\newtheorem{eg}[thm]{Example}

\newcommand{\ben}{\begin{enumerate}}
\newcommand{\een}{\end{enumerate}}
\newcommand{\bpm}{\begin{pmatrix}}
\newcommand{\epm}{\end{pmatrix}}
\newcommand{\bca}{\begin{cases}}
\newcommand{\eca}{\end{cases}}
\newcommand{\beq}{\begin{equation}}
\newcommand{\eeq}{\end{equation}}
\newcommand{\beg}{\begin{eg}}
\newcommand{\eeg}{\end{eg}}
\newcommand{\btm}{\begin{tm}}
\newcommand{\etm}{\end{tm}}
\newcommand{\bpro}{\begin{pro}}
\newcommand{\epro}{\end{pro}}
\newcommand{\brem}{\begin{rem}}
\newcommand{\erem}{\end{rem}}
\newcommand{\bcor}{\begin{cor}}
\newcommand{\ecor}{\end{cor}}

\def\d{\delta}
\def\p{\partial}
\def\b{\beta}
\def\l{\lambda}
\def\a{\alpha}
\def\L{\Lambda}
\def\ti{\tilde}
\def\sdp{\ltimes}

\def\unon{\frac{\U(n)}{\OO(n)}}
\def\runon{\frac{\U(n)\sdp\C^n}{\OO(n)\sdp \R^n}}

\def\ss{\smallskip}
\def\ms{\medskip}
\def\bs{\medskip}
\def\ni{\noindent}
\def\wnnb{with non-degenerate normal bundle}

\begin{document}

\title[Flat Lagrangian immersions]
{Transformations of flat Lagrangian immersions and Egoroff nets}

\author{Chuu-Lian Terng$^1$ \and Erxiao Wang$^2$}\thanks{$^1$ Research supported in part by DMS-0529756 and the UCI Advance Chair fund. $^2$Research supported
in  part by the Postdoctoral Fellowship of MSRI}
\address{University of California, Irvine, CA 92697-3875}
\email{cterng@math.uci.edu}
\address{University of Texas, Austin, TX 78712-0257}
\email{ewang@math.utexas.edu}

\date{}

\begin{abstract}
We associate a natural $\lambda$-family ($\lambda \in \R \setminus \{0\} $) of flat
Lagrangian immersions in $\C^{n}$ with non-degenerate normal
bundle to any given one. We prove that the structure equations for such immersions admit the same 
Lax pair as the first order integrable system associated to 
the symmetric space $\runon$. An interesting observation
is that the family degenerates to an Egoroff net on $\R^n$ when
$\lambda \to 0$.  We construct an action of a rational loop group on such immersions 
by identifying its generators and computing their dressing actions. The action of 
the generator with one simple pole gives the geometric Ribaucour transformation and 
we provide the permutability formula for such transformations. The action of the generator with two poles and the action of a rational loop in the translation subgroup produce new transformations. The corresponding results for flat 
Lagrangian submanifolds in $\C P^{n-1}$ and $\p$-invariant Egoroff nets 
follow nicely via a spherical restriction and 
Hopf fibration. 
\end{abstract}

\maketitle

\section{Introduction}

In recent years techniques of integrable systems have been applied
extensively to study submanifolds and geometric nets. In particular, many curved flat 
systems associated to symmetric space $\frac{U}{K}$ or the gauge 
equivalent $\frac{U}{K}$-systems have shown to be the Gauss-Codazzi-Ricci  equations for 
submanifolds with special geometric properties. To move on, let us briefly go over 
some definitions about geometric nets. 

A vector function $x(u)$, where $x=(x_1, \ldots, x_n)$ are the standard coordinates of $\R^n$ and $u=(u_1, \ldots, u_n)$ are parameters, is called a \emph{net function} when it is a local diffeomorphism. The $u_i$ parameter  lines for all $i$ form a `net' on $\R^n$. When the parameter lines are mutually orthogonal everywhere, $x(u)$ is said to define an 
{\it orthogonal\/} net on $\R^n$ and its inverse function $u(x)$ is called 
an orthogonal (curvilinear) coordinate system on $\R^n$. In particular, the Euclidean metric written in $u$ is diagonal, i.e., $\ud s^2= \sum \|\frac{\p x}{\p u_i} \|^2 \ud u_i^2$. 
This diagonal metric is said to be {\it Egoroff\/} when there exists a \emph{potential function} $\phi(u)$ so that 
$\frac{\p \phi}{\p u_i}=\|\frac{\p x}{\p u_i} \|^2$; then the net $x(u)$ and the  coordinate system $u(x)$ are also called {\it Egoroff \/}. For example, the 
polar coordinate system $u_1=\sqrt{x_1^2+x_2^2}$ and 
$u_2=\arctan (x_2/x_1)$ is orthogonal while $x_1=u_1 \cos u_2$ and 
$x_2=u_1 \sin u_2$ define an orthogonal net on $\R^2$, but they are not 
Egoroff. 

Ribaucour used sphere congruence to construct transformations for
orthogonal coordinate systems in 1872 \cite{Ri}. The
iteration or vectorial extension of Ribaucour transformations was
studied and applied to Egoroff nets on
$\R^n$ by Liu-Ma\~nas in \cite{LM}. Ma\~{n}as, Alonso, and Medina used dressing methods for multicomponent KP hierarchies in \cite{MAM1,MAM2} to construct dressing actions for Egoroff orthogonal nets. Dajczer-Tojeiro generalized
sphere congruence and Ribaucour transformations to submanifolds in 
space-forms in a series of papers \cite{DT2, DT3}, and they found Ribaucour 
transformations for flat Lagrangian submanifolds in $\C^{n}$ and $\C P^n$. 

In this paper we introduce an associated family for any such submanifold which unifies these two geometries, then we generalize both papers \cite{LM} and \cite{DT2} on
Ribaucour transformations not only by discovering new 
types of transformations in closed algebraic formulas, but also by 
describing the group structure of these transformations. The results of 
this paper include:

\ss (1)  We prove that the equation for flat Lagrangian immersions
in $\C^n$ \wnnb\,  is the $\frac{\U(n)\sdp \C^n}{\OO(n)\sdp
\R^n}$-system, where $\U(n)\sdp \C^n$ is the group of unitary rigid
motions of $\C^n$ and $\OO(n)\sdp \R^n$ is the group of orthogonal
rigid motions of $\R^n$.  This equation has a Lax pair, whose frame 
$F(u,\l)$ is of the form $F=\bpm E& X\\ 0&1\epm$.  
We show that for $\l\in \R\setminus \{0\}$, $X(\cdot, \l)$
is a flat Lagrangian immersion  in $\C^n$ \wnnb, while $X(\cdot, 0)$ degenerates 
to an Egoroff orthogonal net on $\R^n$. This $\l$-family of immersions share 
the same induced metric. 

\ss (2) Let $\L^{\tau,\sigma}_{-,m}(n)$ denote the group  of rational
maps $g:S^2=\C\cup \{\infty\}\to \GL(n,\C)$ such that
$g(\bar\l)^*g(\l)=\I$, $g(-\l)^tg(\l)=\I$, and $g(\infty)=\I$. We
find all elements in $\L^{\tau,\sigma}_{-,m}(n)$ with one or two
simple poles (i.e. $g_{\rmi s,\pi_r}$ and $f_{z,\pi}$ respectively in Theorem \ref{an}), and prove that they generate $\L^{\tau,\sigma}_{-,m}(n)$. 

\ss (3) We use dressing action to construct an action of $\L^{\tau,\sigma}_{-,m}(n)$ on flat Lagrangian immersions in $\C P^{n-1}$ and $\p$-invariant flat Egoroff metrics. 

\ss (4) We extend the generators of $\L^{\tau,\sigma}_{-,m}(n)$ to $\L_{-}\GL(n+1,\C)$ and use their dressing actions on the extended frame $F$ to construct transformations of flat Lagrangian immersions in $\C^n$ \wnnb\, and on the space of flat Egoroff metrics on $\R^n$. We then construct dressing actions of rational elements valued in the translation subgroup on  such immersions and metrics. Both produce new type of transformations besides the classical Ribaucour transformations. 

\ss (5) Given a flat Lagrangian immersion $X$ that lies in a hypersphere, it has $\p$-invariant Egoroff metric $\sum_{i=1}^n h_i^2 \ud u_i^2$ and Hopf projection produces flat Lagrangian immersions in $\C P^{n-1}$. We prove that the action in (4) of the extension of $g_{\rmi s, \pi_r}$ on $X$ preserves spherical constraint if and only if $\pi_r(h(0))=0$. However as transformations for such immersions, the method in (3) above produces simpler algorithm. 

\ss (6) As a comparison with the above algebraic dressing actions, we also give the analytic dressing actions, i.e., a first order PDE system that generate a new family of flat Lagrangian immersions from a given one. When the rank of $\pi_r$ is $1$, the action of $g_{\rmi s, \pi_r}$ gives the
geometric Ribaucour transformation constructed  in \cite{DT2}. The
vectorial (or iteration of) Ribaucour transformations in \cite{LM}
for Egoroff nets are given by the action of products of $g_{\rmi s_j, \pi_j}$'s. 

\ss (7) As a bi-product we identify the group structure of Ribaucour
transformations and complex Ribaucour transformations for both flat
Lagrangian immersions in $\C^n$ \wnnb\, and for Egoroff orthogonal 
nets on $\R^n$ in a unified way.

\ss (8) We solve the Cauchy problems for $\unon$-system, flat Egoroff metrics, and flat Lagrangian immersions. 

\ms
The paper is organized as follows. In Section $2$ we review the
geometry of flat Lagrangian submanifolds in $\C^{n}$ and $\C
P^{n-1}$ and introduce the Lax pair, extended frames, and associated
family. In section $3$, we construct generators for $\L^{\tau,\sigma}_{-,m}(n)$, and give explicit formulas for the dressing actions of the generators. In section $4$, we construct an action of $\L^{\tau,\sigma}_{-,m}(n)$ on the space of flat Lagrangian immersions in $\C^n$ that lie in hyperspheres and the space of $\p$-invariant flat Egoroff metrics. In section 5, we compute dressing actions of $g_{\rmi s,\pi_r}$ and of rational elements valued in the translation subgroup on  flat Lagrangian immersions in $\C^n$ \wnnb. In the last section we present the permutability theorem and give explicit formulas for the dressing action of the generator $f_{z,\pi}$. The product of $n$ plane curves in $\C^n$ is the simplest type of flat Lagrangian immersion.  We give an algorithm to construct explicit flat Lagrangian immersions and flat Egoroff metrics by dressing on these product of plane curves.

\section{Geometry of flat Lagrangian submanifolds in $\C^n$ and $\C P^{n-1}$}

We present the geometry of flat Lagrangian submanifolds 
(cf. \cite{DT2, Te2}) in a way to introduce the Lax pair, the
extended frame, and the associated family of them. 

Let $\langle \ , \ \rangle$ and $\omega $ be the standard inner
product and symplectic form on $\C^n = \R^{2n}$ respectively, i.e.,
\[
\langle Z_1, Z_2 \rangle := \re(\bar{Z_1^t} Z_2), \quad
\omega(Z_1,Z_2):= \im(\bar{Z_1^t} Z_2), \quad \forall ~ Z_1, Z_2 \in
\C^n.
\]
Identifying $Z=X+\rmi Y$ in $\C^n$ as $X \choose Y$ in $\R^{2n}$,
the complex structure $J$ on $\R^{2n}$ satisfying
$\omega(\cdot,\cdot) = g(J ~ \cdot, \cdot) $ is:
\[
J \begin{pmatrix} X \\ Y \end{pmatrix} = \begin{pmatrix} -Y \\ X
\end{pmatrix},  \quad \forall ~ X, Y \in \R^n.
\]
Then both $\U(n)$ and its Lie algebra $\uu(n)$ can be identified as
the elements of $\OO(2n)$ and $\oo(2n)$ respectively that commute
with $J$. So
\[
\uu(n)= \left\{
\begin{pmatrix}
P & -Q \\
Q & P
\end{pmatrix}
\in \oo(2n)\ \bigg| \ P\in \oo(n), Q = Q^t \right\}.
\]

A real $n$ dimensional submanifold $M$ in $\C^n$ is
\emph{Lagrangian} if $\omega |_M=0$, or equivalently $J$ maps $TM$
onto the normal bundle $T^{\bot}M$.  The normal bundle $T^{\bot}M$ is {\it
non-degenerate\/} if the dimension of the space of shape operators
at $p$ is equal to codim$(M)$ for all $p\in M$. We call $\U(n)\sdp \R^{2n}$ \emph{unitary rigid motion} group, which preservs the metric and the Lagrangian condition. The following theorem is  known (cf. \cite{DT2, Te2}). Though the proof is elementary, we give it here to set up notations and to be self-contained. 

\begin{tm} \cite{DT2, Te2}
If $f: M \hookrightarrow \R^{2n} \cong \C^n$ is a simply connected flat Lagrangian immersed submanifold  \wnnb, then there exist global line of curvature coordinate 
$u=(u_1,\cdots,u_n)$ on $M$ and a potential function $\phi(u)$ such that the 
fundamental forms of $M$ are:
\begin{equation}\label{av}
\begin{cases}
\I= \sum \frac{\pa \phi}{\pa u_i} \ud u_i^2, &\\
\II= \sum \ud u_i^2 \otimes J(\frac{\pa f}{\pa u_i}).&\\
\end{cases}
\end{equation}
Conversely, given $\phi(u)$ such that $\frac{\pa \phi}{\pa u_i}> 0$ and $\sum \frac{\pa \phi}{\pa u_i} \ud u_i^2$ is flat, the above fundamental forms determine such a flat Lagrangian submanifold uniquely up to unitary rigid motion. 
\end{tm}

\begin{proof}
If $e_1, \ldots, e_n$ is any local orthonormal tangent frame on $M$, then  $W=(e_{1}, \ldots, e_{n}, J e_{1}, \ldots, J e_{n} )$ is a $\U(n)$-valued moving frame and the flat
connection one-form $w:=W^{-1} \ud W$ is $\uu(n)$-valued, i.e.,
\begin{eqnarray}
w_{ij} &=& w_{n+i,n+j},  \label{eqw1}  \\
w_{i,n+j} &=&  w_{j,n+i}, \label{eqw2}
\end{eqnarray}
and $w_{AB}= -w_{BA}$. 
The flatness of $w$,  $dw+w\wedge w=0$, gives the Gauss-Codazzi-Ricci equations  for the immersion $f$, i.e., 
$$\bca 
dw_{ij}=-\sum_{A=1}^{2n} w_{iA}\wedge w_{Aj},\\
dw_{i, n+j}= -\sum_{A=1}^{2n} w_{i A} \wedge w_{A, n+j},\\
dw_{n+i, n+j}=-\sum_{A=1}^{2n} w_{n+i, A}\wedge w_{A, n+j}, 
\eca$$
respectively. Since the induced metric is flat, by \eqref{eqw1} the normal bundle is 
also flat and the Ricci equations become redundant. Recall that a shape operator 
$S_\nu$ on $TM$ is defined by 
$\langle S_\nu (\cdot), \cdot \rangle = \langle \II(\cdot, \cdot), \nu \rangle $ 
for any normal field $\nu$. The flatness of the normal bundle implies that all shape 
operators $\{ S_\nu \} $ commute and thus share eigenspaces. The non-degeneracy 
condition of the normal bundle guarantees 
$n$ independent eigenvectors as the tangent frame. Let $\{e_i\}$ be the normalized 
eigenvectors and $\{\theta_i\}$ the dual $1$-forms. Then $w_{i,n+j}=0$ 
for $i \neq j$, and $w_{i,n+i}= c_i \theta_i$ 
for some smooth functions $c_i$. The non-degeneracy of the normal bundle
implies nonzero principal curvatures $c_i$ along principal 
directions $e_i$. It follows from the Codazzi equations that $\ud (c_i \theta_i)=0$. 
Thus a local line of curvature coordinate $u$ is 
obtained by $c_i \theta_i=\ud u_i$. Introduce $h_i := 1/c_i$, and write 
$\theta_i=h_i \ud u_i$ and their dual $e_i = \frac{1}{h_i}\frac{\pa f}{\pa u_i}$. 
Then the fundamental forms of $M$ are: 
\begin{equation*}
\begin{cases}
\I= \sum h_i^2 \ud u_i^2, &\\
\II= \sum \ud u_i^2 \otimes J(\frac{\pa f}{\pa u_i}).&\\
\end{cases}
\end{equation*}
Moreover, the Levi-Civita connection $1$-form $w_{ij} = \beta_{ij} 
\ud u_i - \beta_{ji} \ud u_j $ with $\beta_{ij}=(h_i)_{u_j}/h_j$ 
for $ i\neq j$ and $\beta_{ii}=0$. In fact the Codazzi equations imply that 
$\beta_{ij}=\beta_{ji}$. This symmetry is clearly equivalent to the existence of a 
potential function $\phi$ making $h_i^2=\frac{\pa \phi}{\pa u_i}$. So the fundamental forms 
are given as in \eqref{av}. When $M$ is simply connected, such a line of curvature 
coordinate $u$ is globally defined. The converse is a direct application of the 
Fundamental Theorem for submanifolds in Euclidean spaces. 
\end{proof}

In the classical literature, $\beta_{ij}$ are called \emph{rotation coefficients}, and a metric taking the special form $\sum \phi_{u_i} \ud u_i^2$ is called an \emph{Egoroff metric}. We call $\b=(\b_{ij})$ the {\it rotation coefficient matrix\/} of the metric $\sum h_i^2 \ud u_i^2$. Then using the eigenframe, we can write simply 
\[
w = W^{-1} \ud W =
\begin{pmatrix}
[\delta, \beta] & -\delta \\
\delta & [\delta, \beta]
\end{pmatrix} , \quad {\rm where\,\, } \delta = \diag(\ud u_1, \ldots, \ud
u_n) . 
\]
The flatness of $w$, i.e., $dw=-w\wedge w$, written in
terms of symmetric $\beta$ is: 
\begin{equation} \label{eqan}
\bca (\beta_{ij})_{u_k}= \beta_{ik}\beta_{kj}, & \textrm{if} \quad  i, j, k  \quad \textrm{are distinct} , \\
(\beta_{ij})_{u_i} + (\beta_{ij})_{u_j} + \sum_k \beta_{ik}
\beta_{jk} =0, & \textrm{if} \quad i \neq j, \eca
\end{equation}
or equivalently, the Levi-Civita connection $1$-form $[\d, \b]$ is flat. These are the {\it
Darboux-Egoroff equations\/} in the classical literature.  It is
also the {\it $\unon$-system\/}, i.e., the first $n$ flows for the
reduced $n$-wave equation in soliton theory (\cite{TU2}).

\ni{\bf Lax pair for the $\unon$-system} \ss

It is known (cf. \cite{Te2}) and easy to check that  $\b$ is a solution of \eqref{eqan} if and only if 
$$w_\lambda = \rmi\lambda \delta + [\delta, \beta]$$
is flat for all  $\lambda\in \C$,
i.e., equation \eqref{eqan} has a {\it Lax pair\/}.

The flatness of $w_\lambda$ implies that there exists  a unique
$E(u,\l)$ satisfying
\begin{equation}\label{eqfra}
E^{-1} \ud E= \rmi \lambda \delta + [\delta, \beta] , \qquad E(0,
\lambda) = I,
\end{equation}
Note that $w_\lambda$ satisfies the following reality condition:
\begin{equation} \label{eqre1}
- w_{\bar\lambda}^{\ast} = w_{\lambda}, \qquad - w_\lambda^{t} = w_{
- \lambda },
\end{equation}
thus its frame $E$ satisfies the $\unon$-reality condition,
\begin{equation} \label{eqre2}
E(u,{\bar\lambda})^{\ast} E(u,{\lambda}) = I, \qquad
E(u,\lambda)^{t} E(u,{ - \lambda }) = I.
\end{equation}
In particular, $E(u,\lambda)\in \U(n)$ for real $\lambda$.

\ms \ni{\bf Extended frame, the associated family and $\runon$-system} \ss

Identify the immersion $f$ and the eigenframe $e_i \in \R^{2n}$ as a complex
column vector in $\C^n$ as before, and set $g=(e_1, \ldots, e_n)$.
Then $g$ is a unitary $n \times n$ matrix and
\begin{equation}\label{aa}
\bca g^{-1} \ud g = \rmi \delta + [\delta, \beta],\\
\ud f= g\d h,\eca
\end{equation}
where $h=(h_1, \ldots, h_n)^t$. Note that \eqref{aa} is solvable for
$g, f$ if and only if the symmetric $\b$ and $h$ satisfy 
\begin{equation}\label{ad}
\bca
[\d, \b]\, {\rm\, is\,\, flat, } \\
(h_i)_{u_j}= \b_{ij} h_j, & i\not= j.
\eca
\end{equation}

The system \eqref{aa} can be also written as
$$k^{-1} \ud k = \bpm \rmi\d + [\d, \b]& \d h\\ 0&0\epm, \quad {\rm where\,\,} k=\bpm g&f\\ 0&1\epm.$$  The following theorem then follows naturally: 

\begin{thm}\label{cc}
Let $\sum_{i=1}^n \frac{\pa \phi}{\pa u_i} \ud u_i^2= \sum_{i=1}^n h_i^2 \ud u_i^2$ 
be a flat Egoroff metric, $\b=(\b_{ij})$ its rotation coefficient matrix, and set 
\begin{equation}\label{bc}
\theta_\l= \bpm \rmi \l \d + [\d, \b]& \d h\\ 0&0\epm.
\end{equation}
  Then: \ben
\item $\theta_\l$ is flat for all $\l\in \C$;
\item There is a unique $F(u,\l)=\bpm E(u,\l)& X(u,\l)\\ 0&1\epm$ solving
$$F^{-1}dF=\theta_\l, \quad F(0,\l)=\I,$$
and $F$ is holomorphic for all $\l\in \C$;
\item $E$ is the frame for $w_\l= \rmi \l \d + [\d,\b]$;  
\item $X_\l= X(\cdot, \l)$ is a flat Lagrangian immersion in $\C^n$ for $\l \in \R$ with
$$\bca \I= \sum_{i=1}^n \frac{\pa \phi}{\pa u_i} \ud u_i^2,\\
\II_\l= \l \sum_{i=1}^n \ud u_i^2 \otimes J(\frac{\p}{\p u_i}),\eca$$ 
and $X(\cdot, 0)\in \R^n \subset \C^n $ is an Egoroff net on $\R^n$. 
\een
\end{thm}

Since the family $X_\l$ shares and only depends on the flat Egoroff metric, we call it 
the {\it associated family\/} for the metric or $\phi$. The frame $F$ in the above 
theorem is called the {\it extended frame\/}, which motivates us to formulate system \eqref{ad} as the $\runon$-system as follows.  

\ss Let $\cg$ be the complexification of $\uu(n)\sdp
\R^{2n}$, i.e.,
\[
\cg := \left\{ \begin{pmatrix}
b & -c & x \\
c & b & y \\
0 & 0 & 0
\end{pmatrix} \Bigg| ~ b=-b^t, ~ c=c^t, ~ b,c \in \gl(n,\C),  ~ x,y \in \C^n  \right\}.
\]
We give two commuting involutions on $\cg$ that gives the symmetric space $\runon$: 
$\tau(A)=\bar{A}$, and $\sigma(A)=TAT$ with 
\[
T:= \begin{pmatrix}
I_n & 0 & 0 \\
0 & -I_n  & 0 \\
0 & 0 & 1
\end{pmatrix} .
\]
It is easy to see that $\sigma\tau= \tau \sigma$,  the fixed point
set of $\tau$ is $\uu(n)\sdp \R^{2n}$, and the fixed point set of
$\sigma$ on $\uu(n)\sdp \R^{2n}$ is $\oo(n)\sdp \R^n$.  So the
corresponding symmetric space is $\frac{\U(n)\sdp
\R^{2n}}{\OO(n)\times \R^n}$. The Cartan decomposition
$\cu=\ck+\cp$ is 
 \begin{align*}
 \ck &=\left\{\bpm b&0&x\\ 0&b& 0\\ 0&0&0\epm\bigg|\, b\in \oo(n), x\in \R^n\right\},\\
 \cp&=\left\{\bpm 0&-c &0\\ c&0&y\\ 0&0&0\epm\bigg|\, c^t=c, \bar c=c, y\in \R^n\right\}.
 \end{align*}
 Then  $\{a_j=e_{n+j,j}-e_{j,n+j}\n 1\leq j\leq n\}$ generates a maximal abelian subalgebra in $\cp$.  The $\frac{U}{K}$-system (cf. \cite{Te1}) is the following PDE
 $$-[a_i, q_{u_j}] + [a_j, q_{u_i}]+[[a_i, q], [a_j, q]]=0, \quad i\not= j,$$ where
\[
q=\begin{pmatrix}
0 & \beta &  0 \\
-\beta & 0 & -h \\
0 & 0 & 0
\end{pmatrix} \in \cp, \quad \textrm{with } \beta_{jj}=0.
\]
Or equivalently,
\[
\theta_{\lambda}= \sum(a_j \lambda+[a_j,q])\ud u_j = \begin{pmatrix}
[\delta, \beta] & \lambda \delta & \delta h \\
-\lambda\delta & [\delta, \beta] & 0 \\
0 & 0 & 0
\end{pmatrix}
\]
is flat for all $\l\in \C$. This is exactly the Lax pair \eqref{bc} for the 
equations of the associated family. So we have proved: 

\btm The  equation for flat Lagrangian immersions in $\C^n$ \wnnb\, is the $\runon$-system with the Lax pair given by \eqref{bc}. \etm

\ms \ni{\bf Flat Lagrangian immersions in $\C^n$ that lie in
$S^{2n-1}$} \ss

\'Elie Cartan proved that an $n$-dimensional flat submanifold can be
locally isometrically  immersed into $S^{2n-1}$, but not in any
lower dimensional spheres. Moreover, the normal bundle of a flat
$n$-dimensional submanifold of $S^{2n-1}$ is automatically flat and
non-degenerate. Thus when a flat Lagrangian immersion $f$ in $\C^n$ lies in a 
hypersphere with center $c_0$ and radius $r$, there always exists the eigenframe $e_i$. Since $f-c_0$ is perpendicular to the tangent plane of $f$, we can write 
$f=c_0+ \sum_{i=1}^n f_i Je_i$ 
for some smooth functions $f_i$. Differentiate it and compare
with $\ud f= \sum_{i=1}^n h_i \ud u_i \otimes e_i$, we get $f_i= h_i$ and $\ud h_i + \sum_{j=1}^n
w_{ji} h_j=0$, i.e.,
\begin{equation}\label{ac}
\ud h + [\d, \b] h=0.
\end{equation}
This implies that $\partial h=0$, where $\p=\sum_{i=1}^n \frac{\p}{\p u_i}.$  
Note also that $\sum h_i^2 \equiv r^2$. In summary, we have  

\bpro\label{ao} Let $\sum_{i=1}^n h_i^2 \ud u_i^2$ be a flat Egoroff
metric, and $\b=(\b_{ij})$ its rotation coefficient matrix.  Then
the following statements are equivalent: \ben
\item $\p h_i=0$, for all $i$, where $\p=\sum_{i=1}^n \frac{\p}{\p u_i}$. 
\item $\ud h+[\d,\b]h=0$.
\item\begin{equation}\label{eqslam}
\begin{cases}
(h_i)_{u_j} = \beta_{ij} h_j ,       \qquad  i \neq j , \\
(h_i)_{u_i} = - \sum \beta_{ij} h_j  .
\end{cases}
\end{equation}
\item $||h||^2=\sum_{i=1}^n h_i^2$ is constant.
\een 
Such flat Egoroff metric will be called \emph{$\p$-invariant or spherical}.
\epro

In the following Theorem, we give an  explicit
formula of the associated family for a $\p$-invariant flat Egoroff metric:

\begin{tm} \label{ap}
Suppose $\sum h_i^2(u) \ud u_i^2$ is a $\p$-invariant flat Egoroff metric, 
$\b$ its rotation coefficient matrix, and $E$ the frame for $\rmi \l\d + [\d,\b]$. Then: 
\ben
\item The associated family in $\C^n$ for the metric is ($\l \in \R \setminus \{0\}$): 
\begin{equation}\label{eqsph}
X(u,\lambda) = -\rmi \lambda^{-1} \left( E(u,\lambda) h(u) - h(0)
\right).
\end{equation}
\item $X(u,\lambda)$ lies in a hyper-sphere
centered at $-\rmi \lambda^{-1} h(0)$ with radius $ \|h(0)\| /
|\lambda| $ for $\l\in\R\setminus \{0\}$.
\item $E(u,0)h(u)= h(0)$, or equivalently, $h(u)= E(u,0)^{-1}h(0)$.
\item $\lim_{\l\to 0} X(u,\l) \in \R^n$ exists and is equal to $-\rmi \frac{\p E}{\p \l}(u,0) h(u)$, which gives a $\p$-invariant Egoroff net on $\R^n$. 
\een
\end{tm}

\begin{proof}
We can check (1) directly that $X$ satisfies 
\beq\label{eqpos} 
\ud X = E\d h , \qquad X(0,\l)=0.
\eeq 
Since $E$ is unitary for $\l\in\R\setminus \{0\}$ and $||h||^2=\sum h_i^2$ is constant, 
(2) follows easily. Now (3) is true because if we
let $A(u)=E(u,0)$, then
$$\ud (Ah)=(\ud A)h + A \ud h = A [\d, \b] h+ A (-[\d, \b]h)=0.$$  So $E(u,0) h(u)= h(0)$ is a constant vector. Thus $h(u)=E(u,0)^{-1}h(0)$. Lastly (4) 
follows from L'Hospital's Rule and the $\unon$ reality condition. 
\end{proof}

\begin{rem}[Cauchy problems for flat Lagrangian immersions]
Suppose $d_1, \ldots, d_n$ are nonzero real constants so that $d_i^2\not= d_j^2$ for all $i\not=j$.  Let $V_n$ denote the space of all real symmetric $n\times n$ matrices with zero diagonal entries.  It was proved in \cite{Te1} that there is an open dense subset $\cs_0(\R, V_n)$ of the space of rapidly decaying smooth maps from $\R$ to $V_n$ such that given $\b_0\in \cs_0(\R, V_n)$ there exists a unique smooth solution $\b$ of the $\unon$-system with initial data $\b(d_1t, \ldots, d_nt)= \b_0(t)$. Specially, 
smooth maps with compact support or with $L^1$ norm less than $1$ lie in $\cs_0(\R, V_n)$.  

Given a solution $\beta$ of the $\unon$-system  \eqref{eqan}, the
following linear system is solvable,
\begin{equation}\label{eqlam}
(h_i)_{u_j} = \beta_{ij} h_j ,       \qquad  i \neq j ,
\end{equation}
and the solutions depend on $n$ nonzero smooth functions of one variable
$b_i(t)$ specifying the initial conditions: $h_i(0, \ldots, u_i, \ldots, 0) = b_i(u_i)$ 
(cf. \cite{Te2}). So the set of Egoroff
metrics with $[\d, \b]$ as its Levi-Civita connections is
parametrized by $n$ positive functions of one variable. In fact, the larger system 
\eqref{eqslam} is still solvable and
the solutions depend only on $h(0)$.  Hence the set of 
$\p$-invariant Egoroff metrics with $[\d, \b]$ as its Levi-Civita 
connection is of finite dimension $n$. They are given by part (3) of Theorem
\ref{ap}. 

Finally Theorem \ref{cc} provides the algorithm to produce the associated family $X(\cdot, \l)$ from $\b$ and $h$. 
\end{rem}

\ms \ni{\bf Flat Lagrangian submanifolds in $\C P^{n-1}$} \ss

If $f$ is a flat Lagrangian immersion in $\C^n$ that lies in a
$S^{2n-1}$ with center $0$, then $f$ is invariant under the $S^1$-action, 
where $S^1$ acts on $S^{2n-1}\subset \C^n$ by $e^{\rmi t}\ast z= e^{\rmi t}z$. 
To see this, we change coordinates linearly from $u$ to $(t_1, \ldots, t_n)$ 
such that $\frac{\p}{\p t_1}= \p = \sum \frac{\p}{\p u_i}$.  Then 
$$E^{-1}\frac{\ud E}{dt_1}= (\rmi \d + [\d, \b])(\frac{\p}{\p t_1})= \rmi I$$
and $f(t_1, t_2, \ldots, t_n)= e^{\rmi t_1} f(0, t_2, \ldots, t_n)$, where $E$ is the frame for $\rmi \d+[\d,\b]$. 

Let $p:S^{2n-1}\to \C P^{n-1}$ denote the Hopf fibration, whose
fibers are orbits of the $S^1$-action on $S^{2n-1}$. Then $N$ is a
flat Lagrangian submanifold of $\C P^{n-1}$ if and only if
$p^{-1}(N)$ is a flat Lagrangian submanifold in $\C^n$ that lies in
$S^{2n-1}$ (\cite{DT2, Te2}). Hence any flat Lagrangian submanifold of $\C P^{n-1}$
comes from a flat Lagrangian submanifold in $\C^n$ that lies in
$S^{2n-1}$.

The explicit assoicated family \eqref{eqsph} of flat Lagrangian submanifolds 
that lie in hyperspheres produces an associated family of $(n-1)$-dimensional flat Lagrangian submanifolds in $\C P^{n-1}$ for $\l \in \R \setminus \{0\}$. 

\section{Dressing actions and generators of rational loop groups}
We construct generators for the group of rational maps $g:\C P^1 \to GL(n,\C)$ that satisfies the reality conditions $g(\bar\l)^*g(\l)=g(-\l)^t g(\l)=\I$ with $g(\infty)=I$, and review the dressing actions of these generators. 

Let $\co_{\infty}$ be an open disk near $\infty$ in $\C P^1 = \C \cup
\{\infty\} $, $\Lambda_{+}(n)$, $\Lambda_{-}(n) $, and $\Lambda(n)$
the groups of holomorphic maps
 from $\C$,  $\co_{\infty}$ and  $\co_{\infty} \cap
\C$ to $\GL(n,\C)$ respectively, and we require $g(\infty)=I$ for
any $g$ in $\Lambda_{-}(n)$. The Birkhoff Factorization Theorem
implies that there exists an open dense subset of $ \Lambda(n) $
such that any $g$ in it can be factored uniquely as $ g = g_{-}
g_{+} = f_{+} f_{-} $ with $ g_{\pm}, f_{\pm} \in \Lambda_{\pm}(n)$.
This open dense subset is called the {\it big cell\/} of $\L(n)$.

The \emph{dressing action} of $\L_-(n)$ on $\L_+(n)$ is defined as
follows: Given $g_\pm \in \L_\pm(n)$,  $g_-\ast g_+= f_+$, where
$f_+$ is the $\L_+(n)$ factor of the factorization $g_-g_+= f_- f_+$
with $f_\pm\in \L_\pm(n)$. Note that this action is local, i.e., it
is defined only when $g_-g_+$ lies in the big cell of $\L(n)$.
However, when $g_-\in \L_-(n)$ is rational, the factorization
$g_-g_+=f_+f_-$ can be computed explicitly as follows: \ben
\item $f_-$ must have the same poles as $g_-$,
\item  the residues of $g_-g_+ f_-^{-1}$ must be zero at poles of $g_-$,
and this leads to a formula for $f_-$, hence we get a formula for $f_+$.
\een The simplest kind of rational element in $\L_-(n)$ is as
follows: Let $\a_i\in \C$, and $\pi$ a projection matrix, i.e.,
$\pi^2=\pi$. Set
\begin{equation} \label{eqsimple}
g_{\alpha_{1}, \alpha_{2}, \pi}(\lambda) := \pi +
\frac{\lambda-\alpha_{2}}{\lambda-\alpha_{1}} (\I-\pi).
\end{equation}
Then $g_{\a_1, \a_2, \pi}\in \L_-(n)$.

First we recall a Theorem in \cite{TU2} that give an explicit
formula for the dressing action of $g_{\a_1, \a_2, \pi}$ on
$\L_+(n)$:

\begin{thm}[ \cite{TU2} ] \label{lmdr}
Let $g_{\alpha_{1}, \alpha_{2}, \pi} $ be as in \eqref{eqsimple},
and $V_{1}$ and $V_{2}$ denote the image of $\pi$ and $I - \pi $
respectively. Assume $f \in \Lambda_{+}(n)$, and
\begin{equation} \label{eqdr}
 f ( \alpha_{1} )^{-1} (V_{1}) \cap f ( \alpha_{2} )^{-1} (V_{2}) = \{0\}.
\end{equation}
Set $\tilde{\pi}$ is the projection onto $f ( \alpha_{1} )^{-1}
(V_{1})$ with respect to
\[
\C^{n} = f( \alpha_{1} )^{-1} (V_{1}) \oplus f ( \alpha_{2} )^{-1}
(V_{2}).
\]
Then \ben
\item $\ti f= g_{\a_1, \a_2, \pi} f g^{-1}_{\a_1, \a_2, \ti\pi}$ is in $\L_+(n)$, i.e.,
$g_{\alpha_{1}, \alpha_{2}, \pi} f = \tilde{f} g_{\alpha_{1},
\alpha_{2}, \tilde{\pi}}$,
\item $\ti f= g_{\a_1, \a_2, \pi}\ast f$.
\een
\end{thm}

Since the proof is rather simple and is typical, we give a sketch
here.  First note that
$$g^{-1}_{\a_1, \a_2, \pi}(\l)= \pi + \frac{\l-\a_1}{\l-\a_2}\, (\I-\pi).$$
So $\ti f(\l)$ is holomorphic in $\C\setminus\{\a_1, \a_2\}$ and has
simple poles at $\a_1, \a_2$.  The residues of $\ti f$ at $\l=\a_1,
\a_2$ are
\begin{align*}
&{\rm Res}(\ti f, \a_1)= (\a_1-\a_2)(\I-\pi)f(\a_1) \ti\pi,\\
& {\rm Res}(\ti f, \a_2)= (\a_2-\a_1)\pi f(\a_2) (\I-\ti\pi).
\end{align*}
It follows from the definition of $\ti \pi$ that both residues are
zero.  Hence $\ti f$ is holomorphic in $\C$, i.e., $\ti f$ lies in
$\L_+(n)$.   This finishes the proof.

\ms

Recall that the frame $E(u,\cdot)$ of the $\frac{\U(n)}{\OO(n)}$-system
\eqref{eqan} is in $\L_+(n)$ satisfying the
$\frac{\U(n)}{\OO(n)}$-reality condition \eqref{eqre2}. Moreover, the
set of $f\in \L(n)$ that satisfy the $\frac{\U(n)}{\OO(n)}$-reality
condition is a subgroup of $\L(n)$.   So we need to consider several
subgroups of $\L_\pm(n)$.
 Let $$\tau: A \to (\bar{A}^{t})^{-1}, \quad \sigma: A \to (A^{t})^{-1} $$ denote the two commuting involutions on $\GL(n,\C)$ determining the symmetric space $\frac{\U(n)}{\OO(n)}$.  Consider the following `twisted' loop groups:
\begin{align*}
\L^\tau(n)&=\{f\in \L(n)\n \tau(f(\bar\l))=f(\l)\},\\
\L_\pm^\tau(n) &= \L_\pm(n)\cap \L^\tau(n),\\
\Lambda^{\tau, \sigma}(n) & = \{f \in \Lambda(n) \n
\sigma(f(\lambda)) = f(-\lambda) , \tau(f(\lambda)) = f(\bar{\lambda})\}, \\
\Lambda_{\pm}^{\tau, \sigma}(n) & = \Lambda^{\tau,
\sigma}(n) \cap \Lambda_\pm (n),
\end{align*}
We will add the subscript `m', such as $\Lambda_{-,m}(n)$, to denote the subgroup of 
rational elements or meromorphic maps from $\C P^1$ to 
$\GL(n,\C)$ in $\Lambda_{-}(n)$. The Birkhoff factorization respects the 
reality conditions (cf. \cite{TU2}), i.e.,  If $g=g_+g_-$ with $g_\pm \in \L_\pm(n)$, then 
\ben
\item $g \in\L^\tau(n)$ implies that $g_\pm \in \Lambda_{\pm}^\tau(n)$,
\item  $g\in \L^{\tau,\sigma}(n)$ implies that $g_\pm\in \L_\pm^{\tau,\sigma}(n)$.
\een As a consequence, the dressing action of $\L_-(n)$ on $\L_+(n)$
restricts to the dressing actions of $\L_-^\tau(n)$ on
$\L^\tau_+(n)$ and $\L_-^{\tau,\sigma}(n)$ on
$\L^{\tau,\sigma}_+(n)$.

\begin{thm}\cite{Uh}\label{bb}
Given $z\in \C\setminus \R$ and a Hermitian projection  $\pi$ of
$\C^{n}$, then 
\[
g_{ z, \pi} (\lambda) :=g_{z, \bar z,\pi}= \pi + \frac{\lambda-
\bar{z}}{\lambda-z} \pi^{\perp} ~~ \in \L_{-,m}^\tau(n),
\]
where $\pi^\perp=\I-\pi$. Moreover, such elements generate $\L_{-,m}^{\tau}(n)$.
\end{thm}

It was proved in \cite{TU2} that the reality condition
$f(\bar\l)^*f(\l)=\I$ implies that \eqref{eqdr} holds. So
$g_{z,\pi}\ast f$ is defined for all $f\in \L^{\tau}_+(n)$ and  we
have

 \begin{thm}\label{am}\cite{TU2}
 Let $f\in \L_+^\tau(n)$, $z\in \C\setminus \{0\}$, $\pi$ a Hermitian projection of $\C^n$, and $\ti f= g_{z,\pi}\ast f$ the dressing action of $g_{z,\pi}$ at $f$. Then
  $\ti f= g_{z,\pi} fg_{z, \ti \pi}^{-1}$ lies in $\L_+^\tau(n)$, where $\ti\pi$ is the Hermitian projection of $\C^n$ onto $f(z)^{-1}(\im\pi)$.
 \end{thm}

It is easy to see that (\cite{TU2}) $g_{z,\pi}\in \Lambda_{-}^{\tau,
\sigma}(n)$ if and only if $z$ is pure imaginary and $\bar\pi=\pi$.
Note that if $g\in \L^{\tau,\sigma}_{-,m}(n)$ has two simple poles
in $\C\setminus(\R\cup \rmi \R)$, then it follows from the reality
condition that the poles of $g$ must be $z, -\bar z$.  So $g=
g_{-\bar z, \rho}g_{z,\pi}$ for some projections $\rho, \pi$. To
find all such $g$'s that lie in $\Lambda_{-,m}^{\tau,\sigma}(n)$, we
need the permutability formula for $\L_-^\tau(n)$:

\begin{thm} [\cite{TU2} Permutability Theorem] \label{properm} 
Given $g_{z_i, \pi_i}$ in $\Lambda_{-,m}^{\tau}(n)$ 
with $z_i \in \C\setminus \R$ ($i=1,2$) and $z_1 \neq z_2, \bar{z}_2 $, let
$\rho_1$ denote the Hermitian projection of $\C^n$ onto $g_{z_2,
\pi_2}(z_1) (\im \pi_1)$, and $\rho_2$ the Hermitian projection onto
$g_{z_1, \pi_1}(z_2) (\im \pi_2)$. Then
$g_{z_2, \rho_2} g_{z_1, \pi_1} = g_{z_1, \rho_1} g_{z_2, \pi_2}$. 
Moreover, such factorization is unique.
\end{thm}

We are ready to prove the analogue of Theorem \ref{bb} for $\L_{-, m}^{\tau, \sigma}(n)$. 

\begin{thm} \label{thmsigma}
Let  $z \in \C \setminus (\R \cup \rmi \R)$, $\pi$ a Hermitian
projection of $\C^{n}$, and $\rho $ the Hermitian projection onto
$g_{z, \pi }(-\bar{z}) (\im \bar{\pi }) $.  Set
\begin{equation} \label{eq2p}
f_{z,\pi} := g_{-\bar{z}, \rho} g_{z, \pi }.
\end{equation}
Then  $f_{z,\pi}\in \L_{-,m}^{\tau, \sigma}(n)$.
\end{thm}

\begin{proof}
Let $z_2=-\bar{z}_1=-\bar{z}$ and $\pi_2=\bar{\pi}_1=\bar{\pi}$. By
the Permutability Theorem, we have $g_{-\bar z, \rho_2} g_{z,\pi}=
g_{z,\rho_1} g_{-\bar z, \bar\pi}$, where
$$\im(\rho_1)= g_{-\bar z, \bar\pi}(z)(\im\pi), \quad \im(\rho_2)= g_{z,\pi}(-\bar z)(\im\bar\pi).$$
But $\overline{g_{z,\pi}(-\bar z)} = g_{-\bar z, \bar\pi}(z)$, so
$\rho_2= \bar\rho_1=\rho$. In other words, we have 
\beq\label{ah}
f_{z,\pi}:= g_{-\bar{z}, \rho } g_{z, \pi } = g_{z, \bar{\rho} }
g_{-\bar{z}, \bar{\pi} } ,  
\eeq 
which implies easily that $f_{z,\pi}$ satisfies $\sigma$-reality condition in addition.
\end{proof}

\begin{thm}\label{an}
The set of $g_{\rmi s,\pi_r}$'s and $f_{z,\pi}$'s generates
$\L^{\tau,\sigma}_{-,m}(n)$, where $s\in \R\setminus  \{0\}$, $z\in
\C\setminus(\R\cup \rmi \R)$,  $\pi_r, \pi$ are Hermitian
projections, and $\bar\pi_r=\pi_r$.
\end{thm}

\begin{proof}
This theorem can be proved by induction on the total degree of a
rational element similar to the proof  of Theorem \ref{bb} in
\cite{Uh} by Uhlenbeck. The details is as follows:

Given $g\in \L_{-,m}^{\tau,\sigma}(n)$, for any complex number $\a
\in \C \setminus \R$ with $g(\a)\neq 0$, there is a unique $k\geq 0$
such that all entries of the matrix
$$ g':=\left(\frac{\l-\a}{\l-\bar\a}\right)^k g $$
have no pole at $\l=\a$ and $g'(\a)\neq 0$. If $g(\a)= 0$, define
$k=0$. Then $\a$ is said to be a zero of $g$ if $g'(\a)$ has nonzero
kernel. Define the order $K$ of the zero of $\det g'$ at $\a$ to be
\emph{the total order of the zeros at the pair} $(\a,\bar\a)$.
Following \cite{Uh}, the \emph{total degree} of $g$ is the sum of
the total order of the zeros at all pairs.

We call both $g_{\rmi s, \pi_r}$ and $f_{z,\pi}$ simple elements. We
will prove that $g$ can be factored as the product of simple
elements  by induction on the total degree of $g$. If the total
degree of $g$ is zero, then $g$ is the constant $\I$.  Suppose the
total degree of $g$ is $m>0$ and the statement is valid when the
total degree is less than $m$. Then pick any zero $\a$ of $g$, we
have the following three cases:

\ss \ni Case (i).  Ker$(g'(\a))= \C^n$ or $g'(\a)=g(\a)= 0$. 

(1) If $\a=-\rmi s$, then $g(\l) = (\l+\rmi s)h(\l)$ for some
meromorphic function $h$ on $S^2$.  It follows from the
$U(n)$-reality condition that
$$g(\l)= \frac{\l+\rmi s}{\l-\rmi s}\ f(\l)$$
for some rational $f:S^2\to GL(n,\C)$.  Since both $g$ and
$\frac{\l+\rmi s}{\l-\rmi s}$ satisfy the $\unon$-reality condition,
so is $f$.  But the total degree of $f$ is $m-n$.  By the inductive
hypothesis, $f$ can be written as products of simple elements.

(2)
If $\a= \bar z\not\in \R\cup \rmi \R$, then since $\overline{g(-\bar
\l)}= g(\l)$,  ${\rm Ker\/}(g(-z))= \overline{{\rm Ker\/}(g(\bar
z))}= \C^n$.  So $g(\l) = (\l+z)(\l-\bar z) \ti g(\l)$ for some
rational $\ti g$.  But $g(\bar\l)^*g(\l)=\I$ implies that
$g(\l)=\frac{(\l+z)(\l-\bar z)}{(\l-z)(\l+\bar z)} \hat g(\l)$ fo
some rational $\hat g(\l)$.  Because both $g$ and
$\frac{(\l+z)(\l-\bar z)}{(\l-z)(\l+\bar z)}$ satisfy the
$\unon$-reality condition, so is $\hat g$. But the total degree of
$\hat g$ is $m-2n$. By induction hypothesis $\hat g$ is the product
of simple elements. 

\ss \ni Case (ii). $\a= \rmi s$ and Ker$(g'(\rmi s))$ is a proper
linear subspace of $\C^n$ of dimension $ k < n$.  Since
$\overline{g(-\bar\l)}= g(\l)$, $\overline{{\rm Ker\/}(g(\rmi s))} =
{\rm Ker\/}(g(\rmi s))$. So the Hermitian projection $\pi$ onto the
orthogonal complement of Ker$(g'(\rmi s))$ is real, i.e.,
$\bar\pi=\pi$. Set
$$\hat g(\l)= g(\l) \left(\pi + \frac{\l+\rmi s}{\l- \rmi s}\ \pi^\perp\right).$$
The total degree of $\hat g$ is $m-k$. By induction hypothesis,
 $\hat g \in \L^{\tau,\sigma}(n)$ factors as a product of simple
 elements. So is $g$.

 \ss\ni
Case (iii). Suppose $\a\not\in \R\cup \rmi \R$, and Ker$(g'(\a))$ is
a proper linear subspace of $\C^n$ of dimension $ k < n$. Let $\rho$
denote the Hermitian projection onto $\overline{ {\rm Ker\/}(g'(\a))
}$. Then $\rho$ determines the Hermitian projection $\pi$ in
Theorem \ref{thmsigma} to define $f_{\bar\a,\pi}= g_{\bar\a,
\bar\rho} g_{-\a, \bar\pi}$ in $\L^{\tau,\sigma}(n)$. Then $\hat g=
g f_{\bar\a,\pi}$ is in $\L^{\tau,\sigma}(n)$ and has total degree
$m+2k$. Now
$$\hat g'(\a)= g'(\a) f_{\bar\a, \pi}(\a)= g'(\a)
\bar\rho g_{-\a, \bar\pi}( \a)=0,$$
 by definition of $\rho$. So $\hat g$ falls into
 the case (i), subcase (2), which implies that there exists $f\in
\L^{\tau,\sigma}_{-,m}(n)$ such that $\hat g=
\frac{(\l-\a)(\l+\bar\a)}{(\l-\bar\a)(\l+ \a)} f$. Note the total
degree of $f$ is now $m+2k-2n <m$. By induction hypothesis, $f$ (and
$g$) can be written as  a product of simple elements.
\end{proof}

\section{Dressing actions on $\p$-invariant flat Egoroff metrics}
We give explicit algorithm to compute the action of $\L^{\tau,\sigma}_{-,m}(n)$ on the space of $\p$-invariant flat Egoroff metrics on $\R^n$, and on the space of flat Lagrangian immersions in $\C^n$ that lie in a hypersphere.  As a consequence, we also get an action of $\L^{\tau,\sigma}_{-,m}(n)$ on the space of flat Lagrangian immersions of $\C P^{n-1}$. 

The dressing action of $\L^{\tau,\sigma}_{-, m}(n)$ on $\L^{\tau,
\sigma}_+(n)$ induces an action on solutions and their frames of the
$\unon$-system. The  Theorem stated below is a consequence of Theorems
\ref{am}, \ref{an}, and Corollary \ref{corbreather}:

\btm\label{ak} Let $\b$ be a solution of the $\unon$-system, and
$E(u,\l)$ its frame, i.e., $E^{-1}\ud E= \rmi \l\d+[\d, \b]$.    Then:
\ben
\item For the first type generator $g_{\rmi \a,\pi}$ of $\L^{\tau,\sigma}_{-, m}(n)$ in 
Theorem \ref{an}, 
\begin{align*}
g_{\rmi \a,\pi}\ast E &=g_{\rmi \a,\pi} Eg_{\rmi \a, \ti\pi}^{-1},\\
g_{\rmi \a,\pi}\ast \b &= \b - 2\a\ti\pi_\ast,
\end{align*}
where  $\xi_\ast:=\xi-\sum_i \xi_{ii} e_{ii}$ and $\ti\pi(u)$ is the
Hermitian projection onto $E(u,\rmi \a)^{-1}(\im \pi)$.
\item For the second type generator $f_{z,\pi}= g_{-\bar z, \rho} g_{z,\pi}$ of $\L^{\tau,\sigma}_{-, m}(n)$, 
\begin{align*}
f_{z,\pi}\ast E&=f_{z,\pi} E g_{z, \hat\pi}^{-1} g_{-\bar z, \hat\rho}^{-1},\\
f_{z,\pi}\ast \b&= \b+ (z-\bar z) (\hat\pi + \hat \rho)_\ast,
\end{align*}
where $\hat \pi$ is the Hermitian projection onto
$E(u,z)^{-1}(\im\pi)$, $\hat\rho(u)$ is the projection onto
$E_1(u,-\bar z)^{-1}(\im\rho)$, and $E_1= g_{z,\pi}E g_{z, \hat
\pi}^{-1}$. \een \etm

\ss
Now Theorem \ref{ap} helps us compute dressing actions on $\p$-invariant flat Egoroff metrics and on flat Lagrangian submanifolds in $S^{2n-1} \subset \C^n$ without using the extended frame.   

\begin{tm}\label{cb}
Let $\ud s^2=\sum h_i(u)^2 \ud u_i^2$ be a $\p$-invariant flat Egoroff metric, $\b$ its 
rotation coefficient matrix, $E(u,\l)$ the frame for $\rmi \l \d+[\d,\b]$, $c=h(0)$, and $X(u,\l)=-\l^{-1}(E(u,\l)E(u,0)^{-1}c-c)$ the associated family of flat Lagrangian immersions given in Theorem \ref{ap}. Let $\ti E= g_{\rmi \a, \pi}\ast E$ and $\hat E= f_{z,\pi}\ast E$ be as in Theorem \ref{ak}. Then constant vectors $\ti c$ and $\hat c$  give the following new $\p$-invariant flat Egoroff metrics and associated family of flat Lagrangian submanifolds in $S^{2n-1} \subset \C^n$: 
\begin{align*}
\ti h(u) &= \ti E(u, 0) \ti c, \quad
\ti X= -\rmi \l^{-1}(\ti E(u,\l)\ti E(u,0)^{-1} \ti c-\ti c),\\
\hat h(u) &= \hat E(u, 0) \hat c,\quad
\hat X = -\rmi \l^{-1}(\hat E(u,\l) \hat E(u, 0)^{-1} \hat c-\hat c). 
\end{align*}
\end{tm}

\bs
\section{Dressing actions on flat Lagrangian immersions in $\C^n$}

We give explicit formulas for dressing actions of rational elements with one simple pole 
on flat Lagrangian immersions in $\C^n$ and on the potential functions of the flat Egoroff metrics, using the extended frame. 

Let $\pi$ be a Hermitian projection of $\C^n$, and $\pi ' $ denote
the Hermitian projection of $\C^{n+1}=\C^n \oplus \C$ onto $\left(
\begin{smallmatrix} \im \pi \\ 0 \end{smallmatrix} \right)$. Then
\begin{equation}\label{eqext}
  g_{z, \bar z, \pi'} = \begin{pmatrix}
  g_{ z, \pi}    &   0 \\
   0   &  \frac{\lambda -\bar{z}}{ \lambda - z }
\end{pmatrix} \in \Lambda_{-}(n+1) .
\end{equation}
We will first compute formally the dressing action of $g_{z, \bar z, \pi'}$ on the extended frame $F$ of a flat Egoroff metric, forgetting the $\sigma$-reality condition $g(-\l)g^t=\I$. 

\begin{thm} \label{lmfac} 
Let $\ud s^2=\sum h_i^2(u) \ud u_i^2$ be a flat Egoroff metric, $\b=(\b_{ij})$ its rotation coefficient matrix, and $F=\bpm E& X\\ 0&1\epm$ its extended frame.  Given $g_{z,\bar z, \pi'}$ as in \eqref{eqext}, then 
$$g_{z, \bar{z}, \pi ' } F = \begin{pmatrix}
    \tilde{E}   &  \tilde{X}  \\
   0   &  1  \end{pmatrix}
   \begin{pmatrix}
  g_{ z, \tilde{\pi}}    &   \xi \\
   0   &  \frac{\lambda -\bar{z}}{ \lambda - z }
\end{pmatrix}\,\, \in \Lambda_{+}(n+1) \times \Lambda_{-}(n+1),$$
where $\tilde{\pi}(u)$ is the Hermitian projection of $\C^{n}$ onto
$E(u,z)^{-1}(\im\pi)$, 
\begin{align*} &\eta = E(u, \bar z)^{-1}X(u, \bar z),\quad \xi= \frac{\bar{z} - z}{ \lambda - z } \tilde{\pi} \eta, \quad \tilde{E} = g_{ z, \pi} E g_{ z, \tilde{\pi}}^{-1},\\
&\tilde{X}= g_{ \bar{z}, \pi^{\perp} } \left(X -
\frac{\bar{z} - z}{ \lambda - z } E \tilde{\pi} \eta \right).\end{align*}
Moreover,  let $\tilde F= \bpm \ti E&\ti X\\ 0&1\epm$, then $\tilde{F}^{-1} \ud \tilde{F}= \begin{pmatrix}
 \rmi\lambda \delta + [\delta, \tilde{\beta}]    & \delta \tilde{h}  \\
  0    &  0
\end{pmatrix}$,
where
\[
\tilde{\beta} = \beta + \rmi(z-\bar{z}) \tilde{\pi}_{\ast}, \quad
\tilde{h}= h +\rmi (z-\bar{z}) \tilde{\pi} \eta.
\]
\end{thm}

\begin{proof}
We first observe that $ F(u, z)^{-1} (\im \pi ') =  \left(
\begin{smallmatrix} E(u, z)^{-1} (\im \pi ) \\ 0 \end{smallmatrix}
\right) $. Since $\im (I-\pi ') = \left( \begin{smallmatrix} \im
\pi^{\perp}  \\ 0 \end{smallmatrix} \right) \oplus \C \left(
\begin{smallmatrix} 0 \\ 1 \end{smallmatrix} \right)$, we have
\[
  F (u, \bar{z} )^{-1}(\im(I-\pi ' )) = \left( \begin{smallmatrix}
E(u, \bar{z})^{-1} (\im \pi^{\perp} ) \\ 0 \end{smallmatrix} \right)
\oplus \C
\left( \begin{smallmatrix} - E(u, \bar{z})^{-1} X(u, \bar{z}) \\
1 \end{smallmatrix} \right)  .
\]
By \eqref{eqre2}, $E(u, z)^{-1}=E(u, \bar{z})^{\ast}$, which implies
that $E(u, z)^{-1} (\im \pi )$ is Hermitian orthogonal to $E(u,
\bar{z})^{-1} (\im \pi^{\perp} )$. Therefore
\[
F (u, z )^{-1} (\im \pi ') \cap F (u, \bar{z} )^{-1}(\im(I-\pi ' ))
= \{0\} .
\]
By Lemma \ref{lmdr}, we can factor
\[
g_{z, \bar{z}, \pi ' } F = \tilde{F} g_{z, \bar{z}, \hat{\pi} } \in
\Lambda_{+}(n+1) \times \Lambda_{-}(n+1) ,
\]
where $\hat{\pi} $ is the projection onto $F (u, z )^{-1}(\im \pi
')$ with respect to
\[
\C^{n+1} = F (u, z )^{-1} (\im \pi ') \oplus F (u, \bar{z}
)^{-1}(\im(I-\pi ' )).
\]
Notice that $\ti \pi$ is no longer a Hermitian projection.  Use
Theorem \ref{lmdr} to get
 $
\hat{\pi} = \begin{pmatrix}
   \tilde{\pi}   &  \tilde{\pi} \eta  \\
   0   &  0
\end{pmatrix}$ and
$g_{z, \bar{z}, \hat{\pi} }  =  \begin{pmatrix}
  g_{ z, \tilde{\pi}}    &   \xi \\
   0   &  \frac{\lambda -\bar{z}}{ \lambda - z }
\end{pmatrix}$, where $\eta= E(u,\bar z)^{-1}X(u,z)$ and  $\xi = \frac{\bar{z} - z}{ \lambda - z } \tilde{\pi} \eta$.

Formulas for $\tilde{E}$ and $\tilde{X}$ can be computed easily from
$\tilde{F} = g_{z, \bar{z}, \pi ' } F g_{z, \bar{z}, \hat{\pi}
}^{-1} $.

Use $\tilde{F} = g_{z, \bar{z}, \pi ' } F g_{z, \bar{z}, \hat{\pi}
}^{-1}$ to compute
\begin{eqnarray} \label{eqlax}
\tilde{F}^{-1} \ud \tilde{F} &=& \begin{pmatrix}
  g_{ z, \tilde{\pi}}    &   \xi \\
   0   &  \frac{\lambda -\bar{z}}{ \lambda - z }
\end{pmatrix}
\begin{pmatrix}
 \rmi\lambda \delta + [\delta, \beta]    & \delta h  \\
  0    &  0
\end{pmatrix}
\begin{pmatrix}
  g_{ z, \tilde{\pi}}    &   \xi \\
   0   &  \frac{\lambda -\bar{z}}{ \lambda - z }
\end{pmatrix}^{-1}
\nonumber\\
&& { } - \ud \begin{pmatrix}
  g_{ z, \tilde{\pi}}    &   \xi \\
   0   &  \frac{\lambda -\bar{z}}{ \lambda - z }
\end{pmatrix} ~ \cdot ~
\begin{pmatrix}
  g_{ z, \tilde{\pi}}    &   \xi \\
   0   &  \frac{\lambda -\bar{z}}{ \lambda - z }
\end{pmatrix}^{-1}.
\end{eqnarray}
Since $\ti F(u,\cdot)\in \L_+(n+1)$, the LHS is holomorphic in
$\l\in \C$. But the RHS has a simple pole at $\l=\infty$.  So $\ti
F^{-1} d\ti F$ must be a degree $1$ polynomial in $\l$.   Use
$(\lambda-z)^{-1} = \lambda^{-1} + z \lambda^{-2} + z^{2}
\lambda^{-3} + \cdots$
to compute the holomorphic part of the RHS to conclude
\[
\tilde{F}^{-1} \ud \tilde{F} = \begin{pmatrix}
 \rmi\lambda \delta + [\delta, \tilde{\beta}]    & \delta \tilde{h}  \\
  0    &  0
\end{pmatrix} ,
\]
with $\tilde{\beta} = \beta + \rmi(z-\bar{z}) \tilde{\pi}_{\ast}$
and $\tilde{h}= h +\rmi (z-\bar{z}) \tilde{\pi} \eta$.
\end{proof}

We will use $g_{z, \bar z, \pi'} \ast F$ to denote $\ti F$, and use $g_{z, \bar z, \pi'} \ast X$ to denote $\ti X$, ... etc.  Now imposing $\sigma$-reality condition, we have:  

\begin{thm} \label{tmnew1} Use the same notation as in Theorem \ref{lmfac}, and 
assume $\bar\pi=\pi$ and $z=\rmi \a$. Then 
\beq\label{bl} \tilde{X} =g_{\rmi \alpha,-\rmi \alpha, \pi'}\ast X = 
g_{-\rmi \alpha, \pi^{\perp} } \left(X + \frac{2\rmi
\alpha}{ \lambda - \rmi \alpha } ~ E ~ \tilde{\pi} \eta \right)
\eeq
is a new associated family of flat Lagrangian immersions for $\ti \I
=\sum_{i} \tilde{h}_{i}^{2} \ud u_{i}^{2}$, where $\tilde{\pi}$ is the Hermitian projection of $\C^{n}$ onto $E(u,\rmi \alpha)^{-1}(\im\pi)$, $\eta= E(u,\rmi \alpha)^{t} 
X(u,-\rmi \alpha)$, and 
$$\tilde{h}=g_{\rmi \alpha,-\rmi \alpha, \pi'} \ast h = h - 2\alpha \tilde{\pi}\eta .$$ 
Moreover, $\ti \I$ is a new flat Egoroff metric with potential
 $$ \tilde{\phi}= g_{\rmi \alpha,-\rmi \alpha, \pi'} \ast \phi = \phi - 2\alpha  \eta^{t} \, \tilde{\pi} \, \eta.$$
Its rotation coefficient matrix is $ \tilde{\beta }= g_{\rmi \alpha,-\rmi \alpha, \pi'} \ast \beta = \beta - 2\alpha \tilde{\pi}_{\ast}$.
\end{thm}

\begin{proof}
Note that the reality conditions \eqref{eqre2} imply $E(u,-\rmi
\alpha)^{-1} = E(u,\rmi \alpha)^{t}$ and $E(u,\rmi \alpha)$ are
real. Using $z=\rmi \alpha$ in Theorem \ref{lmfac}, the uniqueness
of solution to \eqref{eqpos} implies that $X(u,\rmi \alpha)$ and
$\eta$ are real. Thus the new $\tilde{h}$ and $\tilde{\beta}$ are
also real. Hence formulas for $\ti X, \ti h$ and $\eta$ follow from Theorem
\ref{lmfac}. 

To prove the formula for $\ti \phi$, we need to show that $\ud\ti\phi= \sum_i \ti h_i^2
\ud u_i$, i.e., $\ud \ti\phi= \ti h^t\d h$.  Use the formula for $\ti h$
to get $ \tilde{h}^{t} \delta \tilde{h} = \ud \phi - 4\alpha h^{t}
\delta \tilde{\pi} \eta + 4 \alpha^{2} \eta^{t} \tilde{\pi}^{t}
\delta \tilde{\pi} \eta$.  Next we compute $\ud \ti \phi$. Let $V$ be
an $n \times k$ matrix whose columns form a basis for $\im \pi$ and
$U:=E(u,-\rmi \alpha)^{\ast}V$. Then $\tilde{\pi}=U(U^{\ast} U)^{-1}
U^{\ast}$ and
\[
 \tilde{\phi} =  \phi - 2 \alpha X(u,-\rmi \alpha)^{\ast} V (U^{\ast} U)^{-1} V^{\ast} X(u,-\rmi \alpha).
\]
Since $\ud X(u,-\rmi \alpha)^{\ast}= h^{t} \delta E(u,-\rmi
\alpha)^{\ast}$ and $\ud U = (\alpha \delta - [\delta, \beta] )U$,
\begin{eqnarray*}
\ud \tilde{\phi} &=& \ud \phi - 2 \alpha \left[ 2 h^{t} \delta \tilde{\pi} \eta  + X(u,-\rmi \alpha)^{\ast} V \ud(U^{\ast} U)^{-1} V^{\ast} X(u,-\rmi \alpha) \right]  \\
                  &=& \ud \phi - 2 \alpha \left\{ 2 h^{t} \delta \tilde{\pi} \eta - 2 \eta^{\ast} \tilde{\pi}^{\ast} (\alpha \delta - [\delta, \beta] ) \tilde{\pi} \eta  \right\}  \\
                  &=& \ud \phi - 2 \alpha ( 2 h^{t} \delta \tilde{\pi} \eta -2\alpha \eta^{\ast} \tilde{\pi}^{\ast} \delta \tilde{\pi} \eta  )  \\
                  &=& \ud \phi - 4\alpha h^{t} \delta \tilde{\pi} \eta + 4 \alpha^{2} \eta^{t} \tilde{\pi}^{t} \delta \tilde{\pi} \eta  = \tilde{h}^{t} \delta \tilde{h}.
\end{eqnarray*}
We have used $\ud (U^{\ast} U)^{-1} =-(U^{\ast} U)^{-1} \ud
(U^{\ast} U) \, (U^{\ast} U)^{-1}$ in  the second equality above,
and in the third $\eta^{\ast} \tilde{\pi}^{\ast} [\delta,\beta]
\tilde{\pi} \eta = 0$ since $\beta^{\ast}=\beta$.
\end{proof}

\begin{rem}
Our theorem guarantees $\ti X(u,\l)$ in $\ti F$ is holomorphic for $\l\in \C$, though the 
formula for $\ti X$ given by \eqref{bl} seems to have poles at $\l=\pm\rmi \a$. However, 
the residue of \eqref{bl} at $\l=-\rmi \a$ is
\beq\label{bm} -2\rmi \a \pi (X_{-\rmi\a} -E_{-\rmi \a} \ti \pi
\eta) = -2\rmi \a \pi (X_{-\rmi\a} -E_{-\rmi \a} \ti \pi E_{\rmi
\a}^t X_{-\rmi \a}). \eeq 
Substitute $\ti \pi= U(U^*U)^{-1}U$ into \eqref{bm}, and we see that the residue of 
$X_\l$ at $\l=-\rmi \a$ is zero, hence $X_\l$ is holomorphic at $\l=-\rmi \a$. Similar
computation shows that $X_\l$ is holomorphic at $\l=\rmi \a$.
\end{rem} 

By Theorems \ref{ao} and \ref{ap}, a flat Lagrangian immersion in $\C^n$ with
induced metric $\sum h_i^2 \ud u_i^2$ lies in a hypersphere if and only
if $\sum h_j^2$ is constant. In the next Theorem we give a necessary and sufficient 
condition on $\pi$ so that $g_{\rmi \a,-\rmi \a, \pi'} \ast X$ again lies in a hypersphere. 

\begin{tm} \label{tmnew2}
Let $\ud s^2= \sum h_i(u)^2 \ud u_i^2$ be a $\p$-invariant flat Egoroff metric, $\b$ 
its rotation coefficient matrix, and $F=\bpm E& X\\ 0&1\epm$ its extended frame. 
Let $\ti X= g_{\rmi \a,-\rmi \a, \pi'}\ast X$, $\ti h$, $\ti \phi$, $\ti\pi$, and $\eta$ be as in
Theorem \ref{tmnew1}. Then $\ti X$ lie in hyperspheres if and only
if $ \im \pi \perp h(0)$ (or $\pi(h(0))=0$). 

Moreover,  if $ \im \pi \perp h(0)$, then
$\tilde{X}(\cdot,\lambda)$ is contained in the same hypersphere as $X(\cdot,\lambda)$, 
$\alpha \tilde{\pi} \eta = \tilde{\pi} h$, and
\begin{eqnarray}
\tilde{X} &=& g_{-\rmi \alpha, \pi^{\perp} } \left(X + \frac{2\rmi}{ \lambda - \rmi \alpha } ~ E ~ \tilde{\pi} h \right),\label{ca}\\
\tilde{h} &=& h - 2 \tilde{\pi}h, \quad  \tilde{\phi} ~=~ \phi -
\frac{2}{\alpha}  h^{t} \, \tilde{\pi} \, h ,  \quad  \tilde{\beta }
~=~ \beta - 2\alpha \tilde{\pi}_{\ast}.\label{bk}
\end{eqnarray}
\end{tm}

\begin{proof}
Let $V$ and $U$ be the same as in the proof of Theorem \ref{tmnew1}.
Then
$$\pi=V(V^{\ast} V)^{-1} V^{\ast}, \qquad \tilde{\pi}=U(U^{\ast}
U)^{-1} U^{\ast}, \qquad U=E(u,-\rmi \alpha)^{\ast}V.$$
The formula \eqref{eqsph} implies that
$$\eta= E(u,-\rmi \alpha)^{-1} X(u,-\rmi \alpha) = \frac{1}{\alpha}
\left( h - E(u,-\rmi \alpha)^{-1} h(0) \right).$$

If $\im \pi \perp h(0)$, i.e., $h(0)^{t} V =  0$, then by the
formula above we have $\alpha \tilde{\pi} \eta = \tilde{\pi} h$. So
the formulas for $\ti X$ and $\ti h$ are simplified as given.
Moreover, since $I-2\ti\pi= \ti\pi^\perp-\ti \pi$ is orthogonal and
$\ti h= (I-2\ti\pi)h$, we have  $\| \tilde{h} \|^{2} = \| h \|^{2}
=$ constant, which implies that  $\ti X$ lies in a hypersphere.
It follows from $h(0)\perp \im\pi$ and $\ti\pi(0)= \pi$ that
$\tilde{h}(0) = h(0)$. So the new submanifold $\tilde{X}(u,\lambda)$
is contained in the same hypersphere as $X(u,\lambda)$ for each
fixed $\lambda \in \R\setminus \{0\}$.

Conversely, if $||\tilde{h}||^2$ is constant, then  $\tilde{h}^{t}
\ud \tilde{h} = 0$.  We want to prove that $h(0)^{t} V=0$. Use
formula $\tilde{h} = h - 2\alpha \tilde{\pi}\eta$ to compute
\[
\tilde{h}^{t}  \ud \tilde{h} = -2 \alpha (h^{t} - 2\alpha \eta^{t}
\tilde{\pi}) \ud (\tilde{\pi} \eta) - 2\alpha h^{t} [\delta,\beta]
\tilde{\pi} \eta.
\]
Here we have used $\ud h + [\delta,\beta] h = 0$ for spherical case.
Now
\[
\begin{split}
\ud (\tilde{\pi} \eta) &= \ud \left( U(U^{\ast} U)^{-1} V X(u,-\rmi \alpha)   \right)  \\
 &= \left(\alpha \delta - [\delta, \beta] \right) \tilde{\pi} \eta + U \ud (U^{\ast} U)^{-1} \, V X_{-\rmi \alpha} +  \tilde{\pi} \delta h     \\
 &=  \left(\alpha \delta - [\delta, \beta] \right) \tilde{\pi} \eta + \tilde{\pi} \delta \tilde{h}  ,
\end{split}
\]
using $\ud U=\left(\alpha \delta - [\delta, \beta] \right) U $.
After some simplification, we have
\[
\tilde{h}^{t}  \ud \tilde{h} = -2 \alpha (h^{t}-\alpha \eta^{t})
\tilde{\pi} \delta \tilde{h} = -2 \alpha h(0)^{t} V(U^{\ast} U)^{-1}
U^{\ast} \delta \tilde{h}
\]
At $u=0$, this differential is $0$ only when $h(0)^{t} V=0$, since
each entry in $U^{\ast} \delta \tilde{h} |_{u=0} = V^{\ast} \delta h
$ is not $0$.
\end{proof}

As a by-product at $\l=0$, Theorem \ref{tmnew1} and Theorem \ref{tmnew2} give 
formulas for the dressing action of $g_{\rmi \a, \pi}$ on Egoroff and $\pa$-invariant 
Egoroff orthogonal nets respectively. Specially, they give a group point of view of the 
vectorial Ribaucour transformations for Egoroff orthogonal nets constructed in 
\cite{LM}. 

\ms \ni{\bf Dressing actions of general rational elements} \ss

Theorem \ref{lmfac} suggests that we should be able to use the Birkhoff factorization to construct the action of rational loops of the form
$\bpm g_{\rmi \a,\pi} & \frac{b}{\l-\rmi \a}\\ 0 & \frac{\l+\rmi \a}{\l-\rmi \a}\epm$.
This leads us to consider the following element in $\L_-(n+1)$:
Given $\a\in \R\setminus \{0\}$ and $b\in \R^n$, let
 \beq\label{bi} k_{\rmi\a,b}(\l)= \bpm I&
\frac{\rmi b}{\l-\rmi \a}\\ 0& 1\epm. \eeq 
Next we compute the
dressing action of $k_{\rmi \a,b}$ on  flat Lagrangian immersions. Let
$F=\bpm E& X\\ 0&1\epm$ be an extended frame for a flat Lagrangian
immersion $X$ in $\C^n$ \wnnb.  We can factor $k_{\rmi\a, b} F=\ti
F\ti k\in \L_+(n+1)\L_-(n+1)$ as follows:
$$k_{\rmi \a, b} F=\ti F\ti k= \bpm E & Y\\ 0&1\epm \bpm I& \frac{\rmi E_{\rmi \a}^{-1} b}{\l-\rmi\a}\\ 0&1\epm, $$
where
$$Y= X+ \frac{\rmi (b- E_\l E_{\rmi \a}^{-1} b)}{\l-\rmi \a}.$$
Note that $Y$ is holomorphic at $\l=-\rmi \a$, so $\ti F\in
\L_+(n+1)$. A direct computation gives
\begin{align*}
\ti F^{-1}\ud \ti F&=\bpm E^{-1} \ud E& E^{-1} \ud Y\\ 0&0\epm = \bpm \rmi \l \d +[\d, \b]& E^{-1} \ud E + \d E_{\rmi \a}^{-1} b\\ 0&0\epm\\
&=\bpm \rmi \l \d +[\d, \b]& \d h + \d E_{\rmi \a}^{-1} b\\ 0&0\epm.
\end{align*}
Since $b$ is real and $E$ satisfies the $\unon$-reality condition,
we have
$$\overline{E_{\rmi \a}^{-1}b}= \overline{E_{-\rmi \a}^t}b= E_{\rmi \a}^{-1}b.$$ So
$ E_{\rmi \a}^{-1}b$ is real and we have proved

\btm\label{bh}
 Let $X(\cdot, \l)$ be the associated family of flat Lagrangian immersions in $\C^n$ \wnnb\, for the flat Egoroff metric $\sum_{i=1}^n h_i(u)^2 \ud u_i^2$, and $F=\bpm E& X\\ 0&1\epm$ its extended frame.  Let $k_{\rmi \a, b}$ be as in \eqref{bi}. Then
$$k_{\rmi \a, b}\ast X:= X+ \frac{\rmi (b-E_\l E_{\rmi \a}^{-1} b)}{\l-\rmi\a}$$
is a family of flat Lagrangian immersions associated to the flat
Egoroff metric
$$k_{\rmi \a, b}\ast h= h +  E_{\rmi \a}^{-1}b.$$
\etm

Set \beq r_{\rmi \a, \pi, b}(\l):=k_{-\rmi \a, -2\a b}g_{\rmi\a,
-\rmi\a, \pi'}= \bpm g_{\rmi \a, \pi} & \frac{-2\rmi \a b}{\l-\rmi
\a}\\ 0& \frac{\l+\rmi \a}{\l-\rmi \a}\epm. \eeq 
Since both $k_{\rmi
\a, b}\ast F$ and $g_{\rmi\a,-\rmi\a, \pi'}\ast F$ come from the
dressing action of $\L_-(n+1)$ on $\L_+(n+1)$, Theorems \ref{tmnew1}
and \ref{bh} give the formula for the action of $r_{\rmi \a, \pi,b}$
on flat Lagrangian immersions and on flat Egoroff metrics.

\ss \ni{\bf Analytic version of dressing actions} \ss

\ms We can also write down a first order compatible PDE system for
$\ti\pi$ and $b$ that give rise to the action of $r_{\rmi \a,\pi,
b}$ on flat Lagrangian immersions.
 Let $\sum_{i=1}^n h_i^2 \ud u_i^2$ be  a flat Egoroff metric, $\b$ its rotation coefficient matrix, and $F$ its extended frame, i.e.,
$$\theta_\l = F^{-1}dF= \bpm \rmi \d\l + [\d,\b]& \d h\\ 0&0\epm.$$
  We want to find  $\ti F=g F \ti g^{-1}$ with $g=r_{\rmi\a, \pi, b}$ and
$$\ti g= \bpm g_{\rmi \a,\ti\pi(u)}& \frac{-2 \rmi \a}{\l-\rmi \a}\, y(u)\\ 0& \frac{\l+ \rmi \a}{\l-\rmi \a}\epm$$ 
such that
$$\ti\theta_\l=\ti F^{-1}\ud \ti F= \bpm \rmi \l\d+ [\d, \ti \b]& \d\ti h\\ 0&0\epm$$ for some $\ti \b$ and $\ti h$.  Instead of factoring it directly, we can use the fact that
$\ti\theta_\l= \ti g \theta_\l \ti g^{-1} - d\ti g \ti g^{-1}$.  So
$$\ti\theta_\l \ti g= \ti g\theta_\l - d\ti g$$
should hold for all $\l$. This will give a first order PDE system
for $\ti \pi$ and $y$ as follow: If we multiply $(\l-\rmi \a)$ to
both sides of the above equation, then we get
\begin{align*}
& (\rmi \l \d +[\d,\ti\b] (\l-\rmi \a + 2\rmi\a \ti\pi^\perp)\\
&\quad = (\l-\rmi \a + 2\rmi \a\ti\pi^\perp)(\rmi \d\l + [\d, \b]) - 2\rmi \a \ud \ti\pi, \qquad {\rm and}\\
&-2\rmi \a(\rmi \l\d +[\d, \ti\b])y +(\l+\rmi\a)\d\ti h\\
&\quad =(\l-\rmi\a + 2\rmi\a \ti\pi^\perp)\d h + 2\rmi \a dy.
\end{align*}
Equate coefficients of $\l$ and the constant term to get

\btm\label{bj} Suppose $\sum_{i=1}^n h_i(u) \ud u_i^2$ is a flat
Egoroff metric, $\b$ its rotation coefficient matrix, and $X$ the
associated family of flat Lagrangian immersions in $\C^n$ \wnnb.
Then the following system is solvable for $\ti \pi$ and $y$:
\begin{equation}\label{bf}
\bca
\ud \ti\pi= [[\d,\b], \ti\pi] + \a[\d, \ti\pi](\I-2\ti\pi), & \ti\pi(0)=\pi,\\
\ud  y= -[\d, \b-2\a\ti\pi]y +\ti \pi \d h -\a\d y, & y(0)= b, \eca
\end{equation}
where $\ti\pi^2=\ti\pi$, $ \ti\pi^*=\ti\pi =\overline{\ti\pi}$, and
$y(u)\in \R^n$. Moreover, set $\ti h= h-2\a y$, then $\ud \ti
s^2=\sum_{j=1}^n \ti h_j^2\ud u_j^2$ is again a flat Egoroff metric
and $\ti \b= \b- 2\a\ti\pi_\ast$ is the rotation coefficient matrix
for $\ud \ti s^2$. \etm

At the first sight, it is not clear that we can write down all solutions of system \eqref{bf}.  But our algebraic formulas for the actions of $g_{\rmi \a, \pi}$ and $k_{\rmi \a, b}$ in fact gives all the solutions \eqref{bf}. The geometric Ribaucour transformation constructed in \cite{DT2} is just a special case of the above theorem for rank $1$ real projections.

\section{Permutability theorem, complex Ribaucour transformation and  examples}
Let us first compute the permutability formulas for dressing actions of two
simple elements as in Theorem \ref{lmfac}, ignoring the $\sigma$-reality condition. For simplicity, we will use $g_{z, \pi} \ast F$ instead of $g_{z, \bar z, \pi'} \ast F$ in this 
section and so is for $X$. 

\begin{thm} \label{thmperm}
We use the same notations as in Theorem \ref{properm} and Theorem
\ref{lmfac}. For $j=1,2$, let $F_j:=g_{z_j,\pi_j} \ast F$, where
$\ast$ is the dressing action computed in Theorem \ref{lmfac}. Then
we have
\begin{align}\label{gpperm}
  & F_{12}:= ( g_{z_2,\rho_2} g_{z_1,\pi_1} ) \ast F =
g_{z_2,\rho_2} \ast F_{1}   \\
& =F_{21}: =( g_{z_1,\rho_1} g_{z_2,\pi_2} ) \ast F = g_{z_1,\rho_1}
\ast F_{2}.
\end{align}
Moreover,
\begin{align*}
&X_{12}= g_{\bar{z}_2,\rho_2^{\perp}} g_{\bar{z}_1,\pi_1^{\perp}}
\left( X-
\frac{\bar{z}_1 - z_1}{ \lambda - z_1 } E \tilde{\pi}_1 \eta_1 - \frac{\bar{z}_2 - z_2}{ \lambda - z_2 } E ~ g_{z_1,\tilde{\pi}_1^{\perp}} ~ \tilde{\rho}_2 ~ \eta_{12} \right) \\
          &=  g_{\bar{z}_1,\rho_1^{\perp}} g_{\bar{z}_2,\pi_2^{\perp}} \left( X-
\frac{\bar{z}_2 - z_2}{ \lambda - z_2 } E \tilde{\pi}_2 \eta_2 - \frac{\bar{z}_1 - z_1}{ \lambda - z_1 } E ~ g_{z_2,\tilde{\pi}_2^{\perp}} ~ \tilde{\rho}_1 ~ \eta_{21} \right),\\
&F_{12}^{-1}dF_{12}= \bpm i\l\d+ [\d, \b_{12}]& \d h_{12}\\ 0&0\epm,
\end{align*}
where
\begin{align*}
h_{12}&= h + \rmi (z_1-\bar{z}_1) \tilde{\pi}_1 \eta_1 + \rmi
(z_2-\bar{z}_2)
                  \tilde{\rho}_2 ~ \eta_{12} \\
          &=  h + \rmi (z_2-\bar{z}_2) \tilde{\pi}_2 \eta_2 + \rmi (z_1-\bar{z}_1)
                  \tilde{\rho}_1 ~ \eta_{21} , \\
\beta_{12}&= \beta + \rmi (z_1-\bar{z}_1) \tilde{\pi}_{1\ast} + \rmi
(z_2-\bar{z}_2)
                   \tilde{\rho}_{2\ast} \\
          &= \beta + \rmi (z_2-\bar{z}_2) \tilde{\pi}_{2\ast} + \rmi (z_1-\bar{z}_1)
                  \tilde{\rho}_{1\ast}.
\end{align*}
Here $\tilde{\pi}_j (u)$ is the Hermitian projection onto
$E(u,z_j)^{-1}(\im \pi_j)$, $\eta_j (u)= E(u,\bar{z}_j)^{-1}
X(u,\bar{z}_j) $, $\tilde{\rho}_j $ the Hermitian projection onto
$g_{z_{k}, \tilde{\pi}_{k} }(z_j)(\im \tilde{\pi}_j)$, and
$\eta_{kj}(u)=g_{\bar{z}_k, \tilde{\pi}_k^{\perp} }(\bar{z}_j) ~
\eta_j + \frac{\bar{z}_k - z_k}{\bar{z}_k-\bar{z}_j} \,
\tilde{\pi}_k \, \eta_k $ for $j=1,2$ and $k=3-j$.
\end{thm}

\begin{proof}
Because the dressing action is a group action, \eqref{gpperm}
certainly holds by Theorem \ref{properm}. So we have
\[
g_{z_2,\rho_2} ~ E_{1} = E_{12} ~ g_{z_2,\tilde{\rho}_2} ,
\]
for $\tilde{\rho}_2$ being the Hermitian projection onto
$E_{1}(u,z_2)^{-1}(\im \rho_2)$. But
\[
\begin{split}
E_{1}(u,z_2)^{-1}(\im \rho_2) &= E_{1}(u,z_2)^{-1} g_{z_1, \pi_1}(z_2) (\im \pi_2) \\
          &= g_{z_{1}, \tilde{\pi}_{1} }(z_2) E(u,z_2)^{-1}(\im \pi_2) \\
          &= g_{z_{1}, \tilde{\pi}_{1} }(z_2)(\im \tilde{\pi}_2) .
\end{split}
\]
We get $\tilde{\rho}_1$ similarly. Now by Theorem \ref{lmfac},
\[
X_{12}=g_{\bar{z}_2,\rho_2^{\perp}} \left(X_1-  \frac{\bar{z}_2 -
z_2}{ \lambda - z_2 } E_1 ~ \tilde{\rho}_2 ~ \eta_{12} \right),
\]
where
\[
\begin{split}
\eta_{12} &= E_1(u,\bar{z}_2)^{-1} X_1(u,\bar{z}_2) \\
               &= g_{z_1, \tilde{\pi}_1} E^{-1} g_{z_1,\pi_1}^{-1} g_{\bar{z}_1,\pi_1^{\perp}} \left(X - \frac{\bar{z}_1 - z_1}{ \lambda - z_1 } E \tilde{\pi}_1 \eta_1 \right)
               \bigg|_{\lambda=\bar{z}_2} \\
               &= g_{\bar{z}_1, \tilde{\pi}_1^{\perp} }(\bar{z}_2) ~ \eta_2 + \frac{\bar{z}_1 - z_1}{\bar{z}_1-\bar{z}_2} \, \tilde{\pi}_1 \, \eta_1 .
\end{split}
\]
The rest follows directly from Theorem \ref{lmfac}.
\end{proof}

When $z_j$ are pure imaginary and $\pi_j$ are real, the
above theorem gives permutability formulas for the 
transformations in Theorem \ref{tmnew1}. When  $z_2=-\bar{z}_1=-\bar{z}$,
$\pi_2=\bar{\pi}_1=\bar{\pi}$, and $\rho_2=\bar{\rho}_1= \rho$, the
above theorem gives the formula for the dressing action of the other type generator 
$f_{z,\pi}= g_{-\bar{z}, \rho} g_{z, \pi }$ of $\L^{\tau,\sigma}_{-,m}(n)$ on flat Lagrangian
immersions in $\C^n$ \wnnb\, and on Egoroff orthogonal nets. We call 
$X\mapsto f_{z,\pi}\ast X$ a {\it complex Ribaucour transformation\/}. 

\begin{cor} \label{corbreather} Let $\ud s^2=\sum_{i=1}^n h_i^2 \ud u_i^2$ be a flat Egoroff metric, $\b$ the rotation coefficient matrix for $\ud s^2$, and $F=\bpm E& X\\ 0& 1\epm$ the extended frame for $\ud s^2$.
Let $\pi$ be a Hermitian projection of $\C^{n}$, $z \in \C
\setminus (\R \cup \rmi \R)$, and $f_{z,\pi}$ given by Proposition
\ref{thmsigma}.  Then
$$f_{z,\pi}\ast F := g_{-\bar z, \rho}\ast( g_{z,\pi}\ast F)=\bpm \hat E& \hat X\\ 0& 1\epm$$ is an extended frame of a new flat Egoroff metric $\hat\ud s^2=\sum_{i=1}^n \hat h_i^2 \ud u_i^2$. Here 
\begin{align*}
   \hat{X} &= g_{-z,\rho^{\perp}} g_{\bar{z},\pi^{\perp}} \left( X-
\frac{\bar{z} - z}{ \lambda - z} E \tilde{\pi} \eta_1 - \frac{\bar{z} - z}{ \lambda+\bar{z}} E ~ g_{z,\tilde{\pi}^{\perp}} ~ \tilde{\rho} ~ \eta_{12} \right) , \\
   \hat{h} &= f_{z,\pi}\ast h = h + \rmi (z-\bar{z}) ( \tilde{\pi} \eta_1 + \tilde{\rho} ~ \eta_{12}) , \\
   \hat{\beta} &= f_{z,\pi}\ast \b= \beta + \rmi (z-\bar{z})(\tilde{\pi} + \tilde{\rho})_{\ast} ,
 \end{align*}
where $\tilde{\pi}(u)$ is the Hermitian projection of $\C^{n}$ onto
$E(u,z)^{-1}(\im \pi)$, $\tilde{\rho}(u)$ the Hermitian projection
of $\C^{n}$ onto $g_{z, \tilde{\pi} }(-\bar{z})
E(u,-\bar{z})^{-1}(\im \bar{\pi})$, $\eta_1 (u)= E(u,\bar{z})^{-1}
X(u,\bar{z}) $, and $\eta_{12}=g_{\bar{z}, \tilde{\pi}^{\perp}
}E^{-1}X |_{\lambda=-z} + \frac{\bar{z} - z}{\bar{z}+z} \,
\tilde{\pi} \, \eta_1 $.
\end{cor}

\ms
\ni {\bf Flat Lagrangian immersions corresponding to soliton solutions}\par
\ss

The simplest examples of flat Lagrangian submanifolds in
$\C^{n}$ with non-degenerate normal bundle come from the ``vacuum''
solution of $\unon$-system, i.e., $\beta \equiv 0$ or $\beta_{ij}
\equiv 0$. The Lax pair for $\b=0$ is $\rmi \l\d$ with frame  $E = \exp(
\sum_{j}\rmi\lambda u_{j} e_{jj}) $. Let $h_j$ be any smooth positive function defined on some interval $I_j$ around $0$ for $1\leq j\leq n$.  Then $\ud s^2= \sum_{j=1}^n h_j(u_j)^2 \ud u_j^2$ is a flat Egoroff metric with rotation coefficient matrix $\b=0$, and the associated family of $\ud s^2$ is a product of plane curves:
\begin{eqnarray*}
X_\l &=& (z_1, \cdots, z_n )^{t}, \quad \textrm{where } z_j (u_j) = \int_{0}^{u_{j}} h_j(t) e^{ \rmi\lambda t} \ud t ,  \\
\ud s^{2} &=& \sum_{j} \phi_{u_j} \ud u_j^{2} , \quad
\textrm{with } \phi = \sum_{j} \int_{0}^{u_{j}} h_j(t)^{2} \ud t .
\end{eqnarray*}
 
If $\ud s^2$ is $\p$-invariant, then $\sum h_j(u_j)^{2} =c$, which implies
that each $h_{j}$ must be a constant $r_{j}$.  So the corresponding
flat Lagrangian submanifolds of $\C^{n}$ contained in $S^{2n-1}$ (or
$\C P^{n-1}$) are really flat tori, and the associated family is $X_\l=(r_1 (e^{\rmi \l u_1}
-1), \ldots, r_n (e^{\rmi \l u_n}-1))^t / (\rmi \l)$. The potential for the 
$\pa$-invariant Egoroff metric is $\phi=\sum r_{j}^{2} u_{j}$. It is interesting to see when $\l \to 0$ how $X$ gives the standard orthogonal net of $\R^{n}$. 
 
 The orbit of the action of $\L^{\tau,\sigma}_{-,m}(n)$ at the vacuum $\b=0$ is the space of soliton solutions for the $\unon$-system.  We give algorithm to compute flat Lagrangian immersions corresponding to these soliton solutions. 
Since $g_{\rmi \a, \pi_r}$'s and $f_{z,\pi}$'s generate $\L^{\tau,\sigma}_{-,m}(n)$, 
we can use Theorems \ref{tmnew1} and \ref{tmnew2} repeatedly to give a
recursive algorithm to construct explicitly flat Lagrangian immersions
and flat Egoroff metrics given by the action of any $g\in
\L_{-,m}^{\tau,\sigma}(n)$ on these vacuum flat Lagrangian immersions.  We can also 
apply Theorem \ref{bh} to $g\ast X$ and $g\ast h$ to get another $n+1$ parameter families of flat Lagrangian immersions and flat Egoroff metrics.  Similarly, we can apply Theorem \ref{cb} to flat tori repeatedly to get explicit formulas for $\p$-invariant flat Egoroff metrics and on flat Lagrangian immersions in $\C^n$ that lie in hyperspheres.

\bibliographystyle{alpha}


\end{document}